\documentclass[a4paper,12pt]{article}
\begin{document}

\title{Metagroups and their smashed twisted wreath products.}
\author{S.V. Ludkowski}
\date{4 September 2018}
\maketitle
\begin{abstract}
In this article nonassociative metagroups are studied. Different
types of smashed products and smashed twisted wreath products are
scrutinized. Extensions of central metagroups are studied.
\footnote{key words and phrases: metagroup; nonassociative;
product; smashed; twisted wreath  }

\end{abstract}
\par Address: Dep. Appl. Mathematics, Moscow State Techn. Univ. MIREA,
\par av. Vernadsky 78, Moscow 11944, Russia; e-mail: sludkowski@mail.ru

\section{Introduction.}
\par Nonassociative algebras compose a great area of
algebra. In nonassociative algebra, in noncommutative geometry,
quantum field theory, there frequently appear binary systems which
are nonassociative generalizations of groups and related with loops,
quasi-groups, Moufang loops, etc. (see
\cite{bruckb,kakkar,pickert,razm,vojtech} and references therein).
It was investigated and proved in the 20-th century that a
nontrivial geometry exists if and only if there exists a
corresponding loop.
\par Very important role in mathematics and quantum field theory
play octonions and generalized Cayley-Dickson algebras
\cite{albmajja99,allcja98,baez,crimm85,dickson,kansol,schaeferb}.
Their structure and identities in them attract great attention. They
are used not only in algebra and noncommutative geometry, but also
in noncommutative analysis and PDEs, particle physics, mathematical
physics, in the theory of Lie groups and algebras and their
generalizations, mathematical analysis, operator theory and their
applications in natural sciences including physics and quantum field
theory (see
\cite{baez,dickson,frenludkfejms18}-\cite{guertzeb,kansol}-\cite{ludkcvee13}
and references therein). \par A multiplicative law of their
canonical bases is nonassociative and leads to a more general notion
of a metagroup instead of a group \cite{crimm85,ludlmla18}.  The
preposition "meta" is used to emphasize that such an algebraic
object has properties milder than a group. It was used in
\cite{krasner} in French in another meaning corresponding in English
to a multigroup. By their axiomatics metagroups satisfy Conditions
$(2.1.1)$-$(2.1.3)$ and rather mild relations $(2.1.9)$. They were
used in \cite{ludlmla18} for investigations of automorphisms and
derivations of nonassociative algebras.
\par In this article nonassociative metagroups
are studied. Necessary preliminary results on metagroups are
described in Section 2. Quotient groups of metagroups are
investigated in Theorem 2.4. Identities in metagroups established in
Lemmas 2.2, 2.3, 2.6 are applied in Sections 3 and 4. \par Different
types of smashed products of metagroups are investigated in Theorems
3.3 and 3.4. Besides them also direct products are considered in
Theorem 3.1. They provide large families of metagroups (see Remark
3.6).
\par In Section 4 smashed twisted wreath products of metagroups and
particularly also of groups are scrutinized. It appears that
generally they provide loops (see Theorem 4.8). If impose additional
conditions they give metagroups (see Theorem 4.10). In Theorem 4.14
and Corollary 4.15 smashed splitting extensions of non-trivial
central metagroups are studied.
\par All main results of this paper are obtained for the first time.
They can be used for further studies of binary systems,
nonassociative algebra cohomologies, structure of nonassociative
algebras, operator theory and spectral theory over Cayley-Dickson
algebras, PDE, noncommutative analysis, noncommutative geometry,
mathematical physics, their applications in the sciences.

\section{Nonassociative metagroups.}
\par To avoid misunderstandings we give necessary definitions.
For short it will be written metagroup instead of nonassociative
metagroup.
\par {\bf 2.1. Definition.}  Let $G$ be a set with a single-valued
binary operation (multiplication)  $G^2\ni (a,b)\mapsto ab \in G$
defined on $G$ satisfying the conditions: \par $(2.1.1)$ for each
$a$ and $b$ in $G$ there is a unique $x\in G$ with $ax=b$ and \par
$(2.1.2)$ a unique $y\in G$ exists satisfying $ya=b$, which are
denoted by $x=a\setminus b=Div_l(a,b)$ and $y=b/a=Div_r(a,b)$
correspondingly,
\par $(2.1.3)$ there exists a neutral (i.e. unit) element $e_G=e\in G$: $~eg=ge=g$
for each $g\in G$.
\par If the set $G$ with the single-valued
multiplication satisfies the conditions $(2.1.1)$ and $(2.1.2)$,
then it is called a quasi-group. If the quasi-group $G$ satisfies
also the condition $(2.1.3)$, then it is called a loop.
\par The set of all elements $h\in G$
commuting and associating with $G$:
\par $(2.1.4)$ $Com (G) := \{ a\in G: \forall b\in G, ~ ab=ba \} $,
\par $(2.1.5)$ $N_l(G) := \{a\in G: \forall b\in G, \forall c\in G, ~ (ab)c=a(bc) \}
$,
\par $(2.1.6)$ $N_m(G) := \{a\in G: \forall b\in G, \forall c\in G, ~ (ba)c=b(ac)
\} $,
\par $(2.1.7)$ $N_r(G) := \{a\in G: \forall b\in G, \forall c\in G, ~ (bc)a=b(ca)
\} $,
\par $(2.1.8)$ $N(G) := N_l(G)\cap N_m(G)\cap N_r(G)$; \par $Z(G) := Com (G)\cap N(G)$
is called the center $Z(G)$ of $G$.
\par We call $G$ a metagroup if a set $G$ possesses a single-valued binary operation
and satisfies the conditions $(2.1.1)$-$(2.1.3)$ and
\par $(2.1.9)$ $(ab)c=t(a,b,c)a(bc)$ \\ for each
$a$, $b$ and $c$ in $G$, where $t(a,b,c)=t_G(a,b,c)\in Z(G)$.
\par If $G$ is a quasi-group fulfilling also the condition
$(2.1.9)$, then it will be called a strict quasi-group.
\par Then the metagroup $G$ will be called a central metagroup,
if it satisfies also the condition:
\par $(2.1.10)$ $ab={\sf t}_2(a,b)ba$ \\ for each $a$ and $b$ in $G$, where
${\sf t}_2(a,b)\in Z(G)$.
\par If $H$ is a submetagroup (or a subloop) of the metagroup $G$
(or the loop $G$) and $gH=Hg$ for each $g\in G$, then $H$ will be
called almost normal. If in addition $(gH)k=g(Hk)$ and $k(gH)=(kg)H$
for each $g$ and $k$ in $G$, then $H$ will be called a normal
submetagroup (or a normal subloop respectively).
\par Henceforward $Inv_l(a)=Div_l(a,e)$ will denote a left inversion,
\par $Inv_r(a)=Div_r(a,e)$ will denote a right inversion.
\par Elements of a metagroup $G$ will be denoted by small letters,
subsets of $G$ will be denoted by capital letters. If $A$ and $B$
are subsets in $G$, then $A-B$ means the difference of them $A-B=\{
a\in A: ~a \notin B \} $. Henceforward, maps and functions on
metagroups are supposed to be single-valued if something other will
not be specified.

\par {\bf 2.2. Lemma.} {\it If $G$ is a metagroup, then
for each $a$ and $b\in G$ the following identities are fulfilled:
\par $(2.2.1)$ $ ~ b\setminus e=(e/b)t(e/b,b,b\setminus e)$;
\par $(2.2.2)$ $(a\setminus e)b=(a\setminus b)t(e/a,a,a\setminus e)/t(e/a,a,a\setminus b)$;
\par $(2.2.3)$ $b(e/a)=(b/a)t(b/a,a,a\setminus e)/t(e/a,a,a\setminus e) $.}
\par {\bf Proof.} Conditions $(2.1.1)$-$(2.1.3)$ imply that
\par $(2.2.4)$ $b(b\setminus a)=a$, $~b\setminus (ba)=a$;
\par $(2.2.5)$ $(a/b)b=a$, $~(ab)/b=a$ \\
for each $a$ and $b$ in $G$. Using Condition $(2.1.9)$ and
Identities $(2.2.4)$ and $(2.2.5)$ we deduce that
\par $e/b=(e/b)(b(b\setminus e))  = (b\setminus e)
/t(e/b,b,b\setminus e)$\\ which leads to $(2.2.1)$.
\par Let $c=a\setminus b$, then from Identities $(2.2.1)$ and $(2.2.4)$ it
follows that \par $(a\setminus e)b=(e/a)t(e/a,a,a\setminus
e)(ac)$\par $= ((e/a)a)(a\setminus b)t(e/a,a,a\setminus
e)/t(e/a,a,a\setminus b)$\\ which provides $(2.2.2)$.
\par Let now $d=b/a$, then Identities $(2.2.1)$ and $(2.2.5)$ imply
that  \par $b(e/a)=(da)(a\setminus e)/t(e/a,a,a\setminus
e)=(b/a)t(b/a,a,a\setminus e)/t(e/a,a,a\setminus e) $ \\
which demonstrates $(2.2.3)$.

\par {\bf 2.3. Lemma.} {\it Assume that $G$ is a metagroup.
Then for every $a$, $a_1$, $a_2$, $a_3$ in $G$ and $p_1$, $p_2$,
$p_3$ in $Z(G)$:
\par $(2.3.1)$ $t(p_1a_1,p_2a_2,p_3a_3)=t(a_1,a_2,a_3)$;
\par $(2.3.2)$ $t(a,a\setminus e,a)t(a\setminus e,a,e/a)=e$.}
\par {\bf Proof.} Since $(a_1a_2)a_3=t(a_1,a_2,a_3)a_1(a_2a_3)$ and
$t(a_1,a_2,a_3)\in Z(G)$ for every $a_1$, $a_2$, $a_3$ in $G$, then
\par $(2.3.3)$ $t(a_1,a_2,a_3)=((a_1a_2)a_3)/(a_1(a_2a_3))$. \par Therefore, for
every $a_1$, $a_2$, $a_3$ in $G$ and $p_1$, $p_2$, $p_3$ in $Z(G)$
we infer that \par $t(p_1a_1,p_2a_2,p_3a_3)=
(((p_1a_1)(p_2a_2))(p_3a_3))/((p_1a_1)((p_2a_2)(p_3a_3)))$\par $=
((p_1p_2p_3)((a_1a_2)a_3))/((p_1p_2p_3)(a_1(a_2a_3)))=
((a_1a_2)a_3)/(a_1(a_2a_3))$, since \par $(2.3.4)$
$b/(pa)=p^{-1}b/a$ and $b/p=p\setminus b=bp^{-1}$ \\ for each $p\in
Z(G)$, $ ~ a$ and $b$ in $G$, because $Z(G)$ is the commutative
group. Thus $t(p_1a_1,p_2a_2,p_3a_3)=t(a_1,a_2,a_3)$.
\par  From Condition
$(2.1.9)$, Lemma 2.2 and Identity $(2.3.1)$ it follows that \par
$t(a,a\setminus e,a)=((a(a\setminus e))a)/(a((a\setminus e)a))=a/[at(e/a,a,a\setminus e)]$\par $=e/t(a\setminus e,a,e/a)$\\
for each $a\in G$ implying $(2.3.2)$.

\par {\bf 2.4. Theorem.} {\it If $G$ is a metagroup and
$Z_0$ is a subgroup in a center $Z(G)$ such that $t(a,b,c)\in Z_0$
for each $a$, $b$ and $c$ in $G$, then its quotient $G/Z_0$ is a
group.}
\par {\bf Proof.}  As
traditionally the notation is used:
\par $(2.4.1)$ $AB=\{ x=ab: ~ a\in A, ~ b \in B \} $,
\par $(2.4.2)$ $Inv_l(A) = \{ x=a\setminus e: ~ a \in A \} $, \par $(2.4.3)$ $Inv_r(A) = \{
x=e/a: ~ a \in A \} $ \\ for subsets $A$ and $B$ in $G$. Then from
Conditions $(2.1.4)$-$(2.1.8)$ it follows that for each $a$, $b$,
$c$ in $G$ the identities take place
\par $((aZ_0)(bZ_0))(cZ_0)=(aZ_0)((bZ_0)(cZ_0))$ and $aZ_0=Z_0a$.
Evidently $eZ_0=Z_0$. In view of Lemmas 2.2 and 2.3 $(aZ_0)\setminus
e= e/(aZ_0)$, consequently, for each $aZ_0\in G/Z_0$ a unique
inverse $(aZ_0)^{-1}$ exists. Thus the quotient $G/Z_0$ of $G$ by
$Z_0$ is a group.

\par {\bf 2.5. Lemma.} {\it Let $G$ be a metagroup, then $Inv_r(G)$
and $Inv_l(G)$ are metagroups.}
\par {\bf Proof.} At first we consider $Inv_r(G)$. Let $a_1$ and
$a_2$ belong to $G$. Then there are unique $e/a_1$ and $e/a_2$,
since the map $Inv_r$ is single-valued (see Definition 2.1). Since
$Inv_r\circ Inv_l(a)=a$ and $Inv_l\circ Inv_r(a)=a$ for each $a\in
G$, then $Inv_r: G\to G$ and $Inv_l: G\to G$ are bijective and
surjective maps. \par We put $\hat{a}_1\circ
\hat{a}_2=(e/a_2)(e/a_1)$ for each $a_1$ and $a_2$ in $G$, where
$\hat{a}_j=Inv_r(a_j)$ for each $j\in \{ 1, 2 \} $. This provides a
single-valued map from $Inv_r(G)\times Inv_r(G)$ into $Inv_r(G)$.
Then for each $a$, $b$, $x$ and $y$ in $G$ equations $\hat{a}\circ
\hat{x}=\hat{b}$ and $\hat{y}\circ \hat{a}=\hat{b}$ are equivalent
to $(e/x)(e/a)=e/b$ and $(e/a)(e/y)=e/b$ respectively. That is
$\hat{x}=(e/b)/(e/a)$ and $\hat{y}=(e/a)\setminus (e/b)$ are unique.
On the other hand, $e/e=e$, $\hat{e}\circ \hat{b}=e/b=\hat{b}\circ
\hat{e}=\hat{b}$ for each $b\in G$.
\par Then we infer that \par $\hat{a}_1\circ (\hat{a_2}\circ \hat{a}_3)=
((e/a_3)(e/a_2))(e/a_1)=$\par
$t_G(e/a_3,e/a_2,e/a_1)(e/a_3)((e/a_2)(e/a_1))=
t_G(\hat{a}_3,\hat{a}_2,\hat{a}_1)(\hat{a}_1\circ \hat{a}_2)\circ
\hat{a}_3$, \\ consequently,
$t_{Inv_r(G)}(\hat{a}_1,\hat{a}_2,\hat{a}_3)=
e/t_G(\hat{a}_3,\hat{a}_2,\hat{a}_1)$. Evidently, $Inv_r(Z(G))=Z(G)$
and $Z(Inv_r(G))=Z(G)$. Thus Conditions $(2.1.1)$-$(2.1.3)$ and
$(2.1.9)$ are satisfied for $Inv_r(G)$.
\par Similarly putting $Inv_l(a_j)=\check{a}_j$ and $\check{a}_1\circ
\check{a}_2=(a_2\setminus e)(a_1\setminus e)$ for each $a_j\in G$
and $j\in \{ 1, 2, 3 \} $ Conditions $(2.1.1)$-$(2.1.3)$ and
$(2.1.9)$ are verified for $Inv_l(G)$.

\par {\bf 2.6. Lemma.} {\it Assume that $G$ is a metagroup and $a\in
G$, $b\in G$, $c\in G$. Then
\par $(2.6.1)$ $e/(ab)=(e/b)(e/a)t(e/a,a,b)/t(e/b,e/a,ab)$ and
\par $(2.6.2)$ $(ab)\setminus e=(b\setminus e)(a\setminus e)
t(ab,b\setminus e,a\setminus e)/t(a,b,b\setminus e)$. \par $(2.6.3)$
$(a/(bc)=((a/c)/b)t(a/(bc),b,c)$, \par $(2.6.4)$ $(bc)\setminus
a=(c\setminus (b\setminus a))/t(b,c,(bc)\setminus a)$.}
\par {\bf Proof.} From $(2.1.9)$ and $(2.2.5)$ we deduce that
\par $((e/b)(e/a))(ab)=t(e/b,e/a,ab)(e/b)((e/a)(ab))=
t(e/b,e/a,ab)/t(e/a,a,b)$ which implies $(2.6.1)$. Then from
$(2.1.9)$ and $(2.2.4)$ we infer that \par $(ab)((b\setminus
e)(a\setminus e))=((ab)(b\setminus e))(a\setminus e)/t(ab,b\setminus
e,a\setminus e)$\par $= t(a,b,b\setminus e)/t(ab,b\setminus
e,a\setminus e)$ which implies $(2.6.2)$.
\par Utilizing $(2.2.4)$ and $(2.1.9)$ we get $b(c((bc)\setminus
a))=a/t(b,c,(bc)\setminus a)$, hence $c((bc)\setminus a)=(b\setminus
a)/t(b,c,(bc)\setminus a)$ implying $(2.6.4)$.
\par Formulas $(2.2.5)$ and $(2.1.9)$ imply that
$((a/(bc))b)c=t(a/(bc),b,c)a$, consequently,
$(a/c)t(a/(bc),b,c)=(a/(bc))b$ and hence
\par $((a/c)/b)t(a/(bc),b,c)=a/(bc)$.

\section{Smashed products and smashed twisted products of metagroups.}
\par {\bf 3.1. Theorem.} {\it Let $G_j$ be a family of
metagroups (see Definition 2.1), where $j\in J$, $J$ is a set. Then
their direct product $G=\prod_{j\in J}G_j$ is a metagroup and \par
$(3.1.1)$ $Z(G)=\prod_{j\in J}Z(G_j)$.}
\par {\bf Proof.} Each element
$a\in G$ is written as $a= \{ a_j: ~ \forall j\in J, ~ a_j\in G_j
\}$. Therefore, a product $ab= \{ c: ~ \forall j\in J, ~ c_j=a_jb_j,
~ a_j\in G_j,  ~ b_j\in G_j \}$ is a single-valued binary operation
on $G$. Then we get that \par $a\setminus b=\{ d: ~ \forall j\in J,
~ d_j=a_j\setminus b_j, ~ a_j\in G_j,  ~ b_j\in G_j \}$ and \par
$a/b= \{ d: ~ \forall j\in J, ~ d_j=a_j/b_j, ~ a_j\in G_j,  ~ b_j\in
G_j \}$. Moreover, \par $e_G=\{ ~ \forall j\in J, ~ e_{G_j} \}$ is a
neutral element in $G$, where $e_{G_j}$ denotes a neutral element in
$G_j$ for each $j\in J$.  Thus conditions $(2.1.1)$-$(2.1.3)$ are
satisfied.
\par From Conditions $(2.1.4)$-$(2.1.7)$ for each $G_j$ in Definition 2.1
we infer that \par $(3.1.2)$ $Com (G) := \{ a\in G: \forall b\in G,
~ ab=ba \} =$\par $ \{ a\in G: ~ a= \{ a_j: \forall j\in J, a_j\in
G_j \}; \forall b\in G, ~ b= \{ b_j: \forall j\in J, b_j\in G_j \} ;
\forall j\in J, ~ a_jb_j=b_ja_j \} =\prod_{j\in J} Com (G_j)$,
\par $(3.1.3)$ $N_l(G) := \{a\in G: ~ \forall b\in G, ~ \forall c\in G, ~ (ab)c=a(bc)
\} = \{a\in G: ~ a= \{ a_j: \forall j\in J, a_j\in G_j \}; ~ \forall
b\in G, ~ b= \{ b_j: \forall j\in J, b_j\in G_j \}; ~ \forall c\in
G, ~ c= \{ c_j: \forall j\in J, c_j\in G_j \}; ~ \forall j\in J,~
(a_jb_j)c_j=a_j(b_jc_j) \}= \prod_{j\in J} N_l(G_j)$
\\ and similarly \par $(3.1.4)$ $N_m(G)=\prod_{j\in J} N_m(G_j)$ and
\par $(3.1.5)$ $N_r(G)=\prod_{j\in J} N_r(G_j)$.
\par This and $(2.1.8)$ imply that
\par $(3.1.6)$ $N(G)=\prod_{j\in J}N(G_j)$. Thus
\par $(3.1.7)$ $Z(G) := Com (G)\cap N(G)=\prod_{j\in J}Z(G_j)$.
\par Let $a$, $b$ and $c$ be in $G$, then
\par $(ab)c=\{ (a_jb_j)c_j: ~ \forall j \in J, ~ a_j\in G_j, b_j\in
G_j, c_j\in G_j \} $\par $= \{ t_{G_j}(a_j,b_j,c_j) a_j(b_jc_j): ~
\forall j \in J, ~ a_j\in G_j, b_j\in G_j, c_j\in G_j \} =
t_{G}(a,b,c) a(bc)$, where
\par $(3.1.8)$ $t_{G}(a,b,c) = \{ t_{G_j}(a_j,b_j,c_j): ~ \forall j \in J, ~
a_j\in G_j, b_j\in G_j, c_j\in G_j \} $. \par Therefore, Formulas
$(3.1.7)$ and $(3.1.8)$ imply that Condition $(2.1.9)$ also is
satisfied. Thus $G$ is a metagroup.

\par {\bf 3.2. Remark.}
\par $(3.2.1)$. Let $A$ and $B$ be two metagroups and let $Z$
be a commutative group such that $Z_m(A)\hookrightarrow Z$,
$Z_m(B)\hookrightarrow Z$, $Z\hookrightarrow Z(A)$ and
$Z\hookrightarrow Z(B)$, where $Z_m(A)$ denotes a minimal subgroup
in $Z(A)$ containing $t_A(a,b,c)$ for every $a$, $b$ and $c$ in $A$.
\par Using direct products it is always possible to extend either
$A$ or $B$ to get such a case. In particular, either $A$ or $B$ may
be a group. On $A\times B$ an equivalence relation $\Xi $ is
considered such that
\par $(3.2.2)$ $(\gamma v,b)\Xi (v,\gamma b)$ and
$(\gamma v,b)\Xi \gamma (v,b)$ and $(\gamma v,b)\Xi (v,b)\gamma $
\\ for every $v$ in $A$, $b$ in $B$ and $\gamma $ in $Z$.
\par $(3.2.3)$. Let $\phi : A\to {\cal A}(B)$ be a
single-valued mapping, where ${\cal A}(B)$ denotes a family of all
bijective surjective single-valued mappings of $B$ onto $B$
subjected to the conditions $(3.2.4)$-$(3.2.7)$ given below. If
$a\in A$ and $b\in B$, then it will be written shortly $b^a$ instead
of $\phi (a)b$, where $\phi (a) : B\to B$. Let also \par $\eta
_{A,B,\phi }: A\times A\times B\to Z$, $\kappa _{A,B,\phi } :
A\times B\times B\to Z$
\par and $\xi _{A,B,\phi } : ((A\times B)/\Xi ) \times ((A\times B)/\Xi
)\to Z$
\\ be single-valued mappings written shortly as $\eta $, $\kappa $
and $\xi $ correspondingly such that
\par $(3.2.4)$ $(b^u)^v=b^{vu}\eta (v,u,b)$, $~ e^u=e$, $~b^e=b$;
\par $(3.2.5)$ $\eta (v,u,\gamma b)=\eta (v,u,b)$;
\par $(3.2.6)$ $(cb)^u=c^ub^u\kappa (u,c,b)$;
\par $(3.2.7)$ $\kappa (u,\gamma c,b)=\kappa (u,c,\gamma b)=\kappa (u,c,b)$
and \par $\kappa (u,\gamma ,b)=\kappa (u,b,\gamma )=e$;
\par $(3.2.8)$ $\xi ((\gamma u,c),(v,b))=
\xi ((u,c),(\gamma v,b))=\xi ((u,c),(v,b))$ and \par $\xi ((\gamma
,e), (v,b))=e$ and $\xi ((u,c),(\gamma ,e))=e$
\\ for every $u$ and $v$ in $A$, $b$, $c$ in $B$,
$\gamma $ in $Z$, where $e$ denotes the neutral element in $Z$ and
in $A$ and $B$.
\par We put
\par $(3.2.9)$ $(a_1,b_1)(a_2,b_2)=(a_1a_2,\xi
((a_1,b_1),(a_2,b_2))b_1b_2^{a_1})$ \\ for each $a_1$, $a_2$ in $A$,
$b_1$ and $b_2$ in $B$. \par The Cartesian product $A\times B$
supplied with such a binary operation $(3.2.9)$ will be denoted by
$A\bigotimes ^{\phi , \eta , \kappa , \xi }B$.
\par Then we put
\par $(3.2.10)$ $(a_1,b_1)\star (a_2,b_2)=(a_1a_2,\xi
((a_1,b_1),(a_2,b_2))b_2^{a_1}b_1)$ \\ for each $a_1$, $a_2$ in $A$,
$b_1$ and $b_2$ in $B$. \par The Cartesian product $A\times B$
supplied with a binary operation $(3.2.10)$ will be denoted by
$A\star ^{\phi , \eta , \kappa , \xi }B$.

\par {\bf 3.3. Theorem.} {\it Let the conditions of Remark 3.2 be
fulfilled. Then the Cartesian product $A\times B$ supplied with a
binary operation $(3.2.9)$ is a metagroup. Moreover there are
embeddings of $A$ and $B$ into $A\bigotimes ^{\phi , \eta , \kappa ,
\xi }B=C_1$ such that $B$ is an almost normal submetagroup in $C_1$.
If in addition $Z_m(C_1)\subseteq Z_m(B)\subseteq Z$, then $B$ is a
normal submetagroup.}
\par {\bf Proof.} From the conditions of Remark 3.2 it follows
that the binary operation $(3.2.9)$ is single-valued.
\par Let $I_1=((a_1,b_1)(a_2,b_2))(a_3,b_3)$ and $I_2=
(a_1,b_1)((a_2,b_2)(a_3,b_3))$, where $a_1$, $a_2$, $a_3$ belong to
$A$, $b_1$, $b_2$, $b_3$ belong to $B$. Then we infer that \par
$I_1= ((a_1a_2)a_3,\xi ((a_1,b_1),(a_2,b_2)) \xi
((a_1a_2,b_1b_2^{a_1}),(a_3,b_3))(b_1b_2^{a_1})b_3^{a_1a_2})$ and
\par $I_2= (a_1(a_2a_3),\xi ((a_1,b_1),
(a_2a_3,b_2b_3^{a_2}))[\xi ((a_2,b_2),(a_3,b_3))]^{a_1}$\par
$b_1(b_2^{a_1}b_3^{a_1a_2})\kappa (a_1,b_2,b_3^{a_2})\eta
(a_1,a_2,b_3)) $. Therefore
\par $(3.3.1)$ $I_1=t((a_1,b_1),(a_2,b_2),(a_3,b_3))I_2$ with
\par $(3.3.2)$ $t((a_1,b_1),(a_2,b_2),(a_3,b_3))=t_{A}
(a_1,a_2,a_3) t_{B}(b_1,b_2^{a_1},b_3^{a_1a_2})$\par $ \xi
((a_1,b_1), (a_2a_3,b_2b_3^{a_2}))[\xi
((a_2,b_2),(a_3,b_3))]^{a_1}\kappa (a_1,b_2,b_3^{a_2})\eta
(a_1,a_2,b_3)$\par $/[\xi ((a_1,b_1),(a_2,b_2)) \xi
((a_1a_2,b_1b_2^{a_1}),(a_3,b_3))]$.
\\ Apparently
$t_{A\bigotimes ^{\phi , \eta , \kappa , \xi }B}
((a_1,b_1),(a_2,b_2),(a_3,b_3))\in Z$ for each $a_j\in A$, $b_j\in
B$, $j\in \{ 1, 2, 3 \} $, where for shortening of a notation
$t_{A\bigotimes ^{\phi , \eta , \kappa , \xi }B}$ is denoted by $t$.
\par If $\gamma \in Z$, then $\gamma ((a_1,b_1)(a_2,b_2))=
(\gamma a_1a_2,\xi ((a_1,b_1),(a_2,b_2)) b_1b_2^{a_1})= (
a_1a_2,b_1b_2^{a_1})\gamma \xi ((a_1,b_1),(a_2,b_2)) =
((a_1,b_1)(a_2,b_2))\gamma $. Hence $\gamma \in Z(A\bigotimes ^{\phi
, \eta , \kappa , \xi }B)$, consequently, $Z\subseteq Z(A\bigotimes
^{\phi , \eta , \kappa , \xi }B)$.
\par Note that $B^u=B$ for each $u\in A$ by the conditions $(3.2.1)$
and $(3.2.3)$.
\par Next we consider the following equation
\par $(3.3.3)$ $(a_1,b_1)(a,b)=(e,e)$, where $a\in A$, $b\in B$.
\par From $(2.1.2)$ and $(3.2.9)$ we deduce that \par $(3.3.4)$
$a_1=e/a$,
\\ consequently, $\xi ((e/a,b_1),(a,b))b_1b^{(e/a)}=e$
and hence from $(2.1.2)$ and $(3.2.8)$ it follows that \par
$(3.3.5)$ $b_1=e/[\xi ((e/a,b^{(e/a)}),(a,b))b^{(e/a)}]$.
\\ Thus $a_1\in A$ and $b_1\in B$ given by $(3.3.4)$ and $(3.3.5)$
provide a unique solution of $(3.3.3)$.
\par Similarly from the following equation
\par $(3.3.6)$ $(a,b)(a_2,b_2)=(e,e)$, where $a\in A$, $b\in B$,
we infer that
\par $(3.3.7)$ $a_2=a\setminus e$, \\ consequently,
$\xi ((a,b),(a\setminus e,b_2))bb_2^a=e$ and hence $b_2^a= [\xi
((a,b),(a\setminus e,b_2))b]\setminus e$. On the other hand,
$(b_2^a)^{e/a}=\eta (e/a,a,b_2)b_2$, consequently, \par $(3.3.8)$
$b_2= (b\setminus e)^{e/a}/\{ [(\xi ((a,b),(a\setminus e,(b\setminus
e)^{e/a}))]^{e/a}\eta (e/a,a,(b\setminus e)^{e/a}) \} $.
\\ Thus Formulas $(3.3.7)$ and $(3.3.8)$ provide a unique solution of
$(3.3.6)$. \par Notice that $(2.1.3)$ is a consequence of $(3.2.8)$
and $(3.2.9)$.
\par Next we put $(a_1,b_1)=(e,e)/(a,b)$ and
$(a_2,b_2)=(a,b)\setminus (e,e)$ and
\par $(3.3.9)$ $(a,b)\setminus (c,d)=((a,b)\setminus (e,e))(c,d)
$\par $t((e,e)/(a,b),(a,b),((a,b)\setminus (e,e))(c,d))/
t((e,e)/(a,b),(a,b),(a,b)\setminus (e,e))$;
\par $(3.3.10)$ $(c,d)/(a,b)=(c,d)((e,e)/(a,b))$\par $
t((e,e)/(a,b),(a,b),(a,b)\setminus (e,e))/
t((c,d)(e/(a,b)),(a,b),(a,b)\setminus (e,e))$ \par and $e_G=(e,e)$,
where $G=A\bigotimes ^{\phi , \eta , \kappa , \xi }B$. Therefore
Properties $(2.1.1)$-$(2.1.3)$ and $(2.1.9)$ are fulfilled for
$A\bigotimes ^{\phi , \eta , \kappa , \xi }B$.
\par Naturally $A$ is embedded into $C_1$ as $ \{ (a,e): ~ a\in A \}
$ and $B$ is embedded into $C_1$ as $ \{ (e,b): ~ b\in B \} $. Let
$a\in A$ and $b_0\in B$, then $(a,b_0)B= \{ (a, \xi
((a,b_0),(e,b))b_0b^a): ~ b\in B \} $ and $B(a,b_0)= \{ (a, \xi
((e,b),(a,b_0))bb_0): ~ b\in B \} $, since $b_0^e=b_0$ by $(3.2.4)$.
From $B^a=B$, $b_0B=B$, $Bb_0=B$, $Z\subset Z(B)$ and $(3.2.3)$,
$(3.2.8)$ it follows that $(a,b_0)B=B(a,b_0)$, where $B^a= \{ b^a: ~
b\in B \} $. Thus $B$ is an almost normal submetagroup in $C_1$ (see
Definition 2.1). If in addition $Z_m(C_1)\subseteq Z_m(B)\subseteq
Z$, then evidently $B$ is a normal submetagroup (see also Condition
$(3.2.2)$), since $t_{C_1}(g,b,h)\in Z_m(C_1)$ and
$t_{C_1}(h,g,b)\in Z_m(C_1)$ for each $g$ and $h$ in $G$, $b\in H$.
\par {\bf 3.4. Theorem.} {\it Suppose that the conditions of Remark 3.2
are satisfied. Then the Cartesian product $A\times B$ supplied with
a binary operation $(3.2.10)$ is a metagroup. Moreover there exist
embeddings of $A$ and $B$ into $A\star ^{\phi , \eta , \kappa , \xi
}B=C_2$ such that $B$ is an almost normal submetagroup in $C_2$. If
additionally $Z_m(C_2)\subseteq Z_m(B)\subseteq Z$, then $B$ is a
normal submetagroup.}
\par {\bf Proof.} The conditions of Remark 3.2 imply that
the binary operation $(3.2.10)$ is single-valued.
\par We consider the following formulas
\par $I_1=((a_1,b_1)\star (a_2,b_2))\star (a_3,b_3)$ and $I_2=
(a_1,b_1)\star ((a_2,b_2)\star (a_3,b_3))$, where $a_1$, $a_2$,
$a_3$ are in $A$, $b_1$, $b_2$, $b_3$ are in $B$. Utilizing
$(3.2.4)$-$(3.2.8)$ and $(3.2.10)$ we get that
\par $I_1= ((a_1a_2)a_3,\xi ((a_1a_2,b_2^{a_1}b_1),(a_3,b_3)) \xi
((a_1,b_1),(a_2,b_2))b_3^{a_1a_2}(b_2^{a_1}b_1))$ and
\par $I_2= (a_1(a_2a_3),\xi ((a_1,b_1),
(a_2a_3,b_3^{a_2}b_2))[\xi ((a_2,b_2),(a_3,b_3))]^{a_1}$\par
$(b_3^{a_1a_2}b_2^{a_1})b_1\kappa (a_1,b_3^{a_2},b_2)\eta
(a_1,a_2,b_3)) $. Therefore
\par $(3.4.1)$ $I_1=t((a_1,b_1),(a_2,b_2),(a_3,b_3))I_2$ with
\par $(3.4.2)$ $t((a_1,b_1),(a_2,b_2),(a_3,b_3))=t_{A}
(a_1,a_2,a_3) \xi ((a_1,b_1),(a_2,b_2))$\par $\xi
((a_1a_2,b_2^{a_1}b_1),(a_3,b_3)) /\{
t_{B}(b_3^{a_1a_2},b_2^{a_1},b_1)$\par $ \xi ((a_1,b_1),
(a_2a_3,b_3^{a_2}b_2))[\xi ((a_2,b_2),(a_3,b_3))]^{a_1}\kappa
(a_1,b_3^{a_2},b_2)\eta (a_1,a_2,b_3)) \} $,
\\ consequently, $t((a_1,b_1),(a_2,b_2),(a_3,b_3))\in Z$
for each $a_j\in A$, $b_j\in B$, $j\in \{ 1, 2, 3 \} $. We denote
$t((a_1,b_1),(a_2,b_2),(a_3,b_3))$ in more details by \\ $t_{A\star
^{\phi , \eta , \kappa , \xi }B} ((a_1,b_1),(a_2,b_2),(a_3,b_3))$
(see Formula $(3.4.2)$). \par Evidently, $(2.1.3)$ is a consequence
of $(3.2.8)$ and $(3.2.10)$.
\par Note that if $\gamma \in Z$, then
\par $\gamma ((a_1,b_1)\star (a_2,b_2))= (\gamma a_1a_2,\xi
((a_1,b_1),(a_2,b_2)) b_2^{a_1}b_1)=$\par $ (
a_1a_2,b_2^{a_1}b_1)\gamma \xi ((a_1,b_1),(a_2,b_2)) =
((a_1,b_1)\star (a_2,b_2))\gamma $.
\\ Therefore $\gamma \in Z(A\star ^{\phi , \eta , \kappa , \xi }B)$,
consequently, $Z\subseteq Z(A\star ^{\phi , \eta , \kappa , \xi
}B)$.
\par Then we seek a solution of  the following equation
\par $(3.4.3)$ $(a_1,b_1)\star (a,b)=(e,e)$, where $a\in A$, $b\in B$.
\par From $(2.1.2)$ and $(3.2.10)$ it follows that \par $(3.4.4)$
$a_1=e/a$,
\\ consequently, $\xi ((e/a,b_1),(a,b))b^{(e/a)}b_1=e$.
Therefore $(2.1.1)$ and $(3.2.8)$ imply that \par $(3.4.5)$
$b_1=[\xi ((e/a,b^{(e/a)}),(a,b))b^{(e/a)}]\setminus e$.
\\ Thus $a_1\in A$ and $b_1\in B$ prescribed by Formulas
$(3.4.4)$ and $(3.4.5)$ provide a unique solution of $(3.4.3)$.
\par Analogously for the following equation
\par $(3.4.6)$ $(a,b)(a_2,b_2)=(e,e)$, where $a\in A$, $b\in B$,
we deduce that
\par $(3.4.7)$ $a_2=a\setminus e$, \\ consequently,
$\xi ((a,b),(a\setminus e,b_2))b_2^ab=e$ and hence $b_2^a= e/[\xi
((a,b),(a\setminus e,b_2))b]$. From $(3.2.4)$ and $(3.2.5)$ it
follows that  $(b_2^a)^{e/a}=\eta (e/a,a,b_2)b_2$, consequently,
\par $(3.4.8)$ $b_2= (e/b)^{e/a}/\{ [(\xi
((a,b),(a\setminus e,(e/b)^{e/a}))]^{e/a}\eta (e/a,a,(e/b)^{e/a}) \}
$.
\\ Thus a unique solution of
$(3.4.6)$ is given by Formulas $(3.4.7)$ and $(3.4.8)$.
\par Then we put $(a_1,b_1)=(e,e)/(a,b)$ and
$(a_2,b_2)=(a,b)\setminus (e,e)$ and
\par $(3.4.9)$ $(a,b)\setminus (c,d)=((a,b)\setminus (e,e))(c,d)
$\par $t((e,e)/(a,b),(a,b),((a,b)\setminus (e,e))(c,d))/
t((e,e)/(a,b),(a,b),(a,b)\setminus (e,e))$;
\par $(3.4.10)$ $(c,d)/(a,b)=(c,d)((e,e)/(a,b))$\par $
t((e,e)/(a,b),(a,b),(a,b)\setminus (e,e))/
t((c,d)(e/(a,b)),(a,b),(a,b)\setminus (e,e))$ \par and $e_G=(e,e)$,
where $G=A\star ^{\phi , \eta , \kappa , \xi }B$. This means that
Properties $(2.1.1)$-$(2.1.3)$ and $(2.1.9)$ are fulfilled for
$A\star ^{\phi , \eta , \kappa , \xi }B$.
\par Evidently there are embeddings of $A$ and $B$ into
$C_2$ as $(A,e)$ and $(e,B)$ respectively. Suppose that $a\in A$ and
$b_0\in B$, then \par $(a,b_0)\star B = \{ (a, \xi ((a,b_0),(e,b))
b^ab_0): ~ b\in B \} $ and \par $B\star (a,b_0) = \{ (a, \xi
((e,b),(a,b_0)) b_0b): ~ b\in B \} $. \\ Therefore $(a,b_0)\star B =
B\star (a,b_0)$ by the conditions $(3.2.3)$ and $(3.2.8)$, since
$B^a=B$ and $Z\subset Z(B)$. Thus $B$ is an almost normal
submetagroup in $C_2$( see Definition 2.1). If additionally
$Z_m(C_2)\subseteq Z_m(B)\subseteq Z$, then apparently $B$ is a
normal submetagroup (see also Condition $(3.2.2)$), since
$t_{C_2}(g,b,h)\in Z_m(C_2)$ and $t_{C_2}(h,g,b)\in Z_m(C_2)$ for
each $g$ and $h$ in $G$, $b\in B$.

\par {\bf 3.5. Definition.} The metagroup
$A\bigotimes ^{\phi , \eta , \kappa , \xi }B$ provided by Theorem
3.3 (or $A\star ^{\phi , \eta , \kappa , \xi }B$ by Theorem 3.4) we
call a smashed product (or a smashed twisted product
correspondingly) of metagroups $A$ and $B$  with smashing factors
$\phi $, $\eta $, $\kappa $ and $\xi $.

\par {\bf 3.6. Remark.} From Theorems 3.1, 3.3 and 3.5 it follows that
taking nontrivial $\eta $, $\kappa $ and $\xi $ and starting even
from groups with nontrivial $Z(G_j)$ or $Z(A)$ it is possible to
construct new metagroups with nontrivial $Z(G)$ and ranges
$t_{G}(G,G,G)$ of $t_{G}$ may be infinite. \par With suitable
smashing factors $\phi $, $\eta $, $\kappa $ and $\xi $ and with
nontrivial metagroups or groups $A$ and $B$ it is easy to get
examples of metagroups in which $e/a\ne a\setminus e$ for an
infinite family of elements $a$ in $A\bigotimes ^{\phi , \eta ,
\kappa , \xi }B$ or in $A\star ^{\phi , \eta , \kappa , \xi }B$.
Evidently smashed products and smashed twisted products (see
Definition 3.5) are nonassociative generalizations of semidirect
products. Combining Theorems 3.3 and 3.5 with Lemmas 2.5 and 2.6
provides other types of smashed products by taking $\hat{b}_1\circ
\hat{b}_2^{a_1}$ instead of $b_1b_2^{a_1}$ or
$\check{b}_2^{a_1}\circ \check{b}_1$ instead of $b_2^{a_1}b_1$ on
the right sides of $(3.2.9)$ and $(3.2.10)$ correspondingly, etc.
\section{Smashed twisted wreath products of metagroups.}
\par {\bf 4.1. Lemma.} {\it Let $D$ be a metagroup and $A$ be a submetagroup
in $D$. Then there exists a subset $V$ in $D$ such that $D$ is a
disjoint union of $vA$, where $v\in V$, that is,
\par $(4.1.1)$ $D=\bigcup_{v\in V} vA$ and
\par $(4.1.2)$ $(\forall v_1\in V, ~ \forall v_2\in V, ~v_1\ne v_2)
\Rightarrow (v_1A\cap v_2A)=\emptyset )$.}
\par {\bf Proof.} The cases $A= \{ e \} $ and $A=D$ are trivial.
Let $A\ne \{ e \} $ and $A\ne D$ and let $Z(D)$ be a center of $D$.
From the conditions $(2.1.4)$-$(2.1.8)$ it follows that $z\in
Z(D)\cap A$ implies $z\in Z(A)$. \par Assume that $b\in D$ and $z\in
Z(D)$ are such that $zbA\cap bA\ne \emptyset $. It is equivalent to
$(\exists s_1\in A, ~ \exists s_2\in A, ~ zbs_1=bs_2)$. From Formula
$(2.2.4)$ it follows that $(zbs_1=bs_2)\Leftrightarrow
(zs_1=s_2)\Leftrightarrow (z=s_2/s_1\in A)$, because $z\in Z(D)$.
Thus \par $(4.1.3)$ $(\exists b\in D, ~\exists z\in Z(D), ~ zbA\cap
bA=\emptyset )\Leftrightarrow (\exists b\in D, ~ \exists z\in
Z(D)-A)$. \par Suppose now that $b_1\in D$, $b_2\in D$ and $b_1A\cap
b_2A\ne \emptyset $. This is equivalent to $(\exists s_1\in A, ~
\exists s_2\in A, ~ b_1s_1=b_2s_2 )$. By the identity $(2.2.5)$ the
latter is equivalent to $b_1=(b_2s_2)/s_1$. On the other hand,
\par $(b_2s_2)/s_1=(b_2s_2)(e/s_1)t(e/s_1,s_1,s_1\setminus
e)/t((b_2s_2)/s_1,s_1,s_1\setminus e)$\par
$=b_2(s_2(e/s_1))t(b_2,s_2,e/s_1)t(e/s_1,s_1,s_1\setminus
e)/t((b_2s_2)/s_1,s_1,s_1\setminus e)$ \\ by $(2.1.9)$, $(2.2.3)$
and $(2.2.5)$. Together with $(4.1.3)$ this gives the equivalence:
\par $(4.1.4)$ $(\exists b_1\in D, ~ \exists b_2\in D, ~ b_1A\cap
b_2A\ne \emptyset ) \Leftrightarrow $\par $(\exists b_1\in D, ~
\exists b_2\in D, ~ \exists s\in A, \exists z\in Z(D)-A, ~ b_1=zb_2s
)$.
\par Let $\Upsilon $ be a family of subsets $K$ in $D$ such that
$k_1A\cap k_2A=\emptyset $ for each $k_1\ne k_2$ in $K$. Let
$\Upsilon $ be directed by inclusion. Then $\Upsilon \ne \emptyset
$, since $A\subset D$ and $A\ne D$. Therefore, from $(4.1.3)$ and
$(4.1.4)$ and the Kuratowski-Zorn lemma (see \cite{kunenb}) the
assertion of this lemma follows, since a maximal element $V$ in
$\Upsilon $ gives $(4.1.1)$ and $(4.1.2)$.
\par {\bf 4.2. Definition.} A set $V$ from Lemma 4.1 is called a
right transversal (or complete set of right coset representatives)
of $A$ in $D$.
\par The following corollary is an immediate consequence of Lemma
4.1.
\par {\bf 4.3. Corollary.} {\it Let $D$ be a metagroup, $A$ be a
submetagroup in $D$ and $V$ a right transversal of $A$ in $D$. Then
\par $(4.3.1)$ $\forall a\in D, ~ \exists _1 s\in A, ~ \exists _1 b\in
V, ~ a=sb$ for a given triple $(A,D,V)$.}
\par {\bf 4.4. Remark.} We denote $b$ in the decomposition $(4.3.1)$
by $b=\tau (a)=a^{\tau }$ and $s=\psi (a)=a^{\psi }$, where $\tau $
and $\psi $ is a shortened notation of $\tau _{A,D,V}$ and $\psi
_{A,D,V}$ respectively. That is, there are single-valued maps \par
$(4.4.1)$ $\tau : D\to V$ and $\psi : D\to A$. \par These maps are
idempotent $\tau (\tau (a))=\tau (a)$ and $\psi (\psi (a))=\psi (a)$
for each $a\in D$.
\par $(4.4.2)$. If $b=a^{\tau }$, then we denote $e/b$ by $a^{e/\tau
}$ and $b\setminus e$ by $a^{\tau \setminus e}$.
\par According to $(2.1.2)$ $s=a/b$ hence $a^{\psi }=a/a^{\tau }$.
From Formula $(2.2.3)$ it follows that $a/b=a(e/b)t(e/b,b,b\setminus
e)/t(a/b,b,b\setminus e)$, consequently, by Lemma 2.3
\par $(4.4.3)$ $s=aa^{e/\tau }t(a^{e/\tau },a^{\tau },a^{\tau \setminus
e})/t(aa^{e/\tau },a^{\tau },a^{\tau \setminus e})$.
\par Notice that the metagroup need not be power-associative.
Then $e/s$ and $s\setminus e$ can be calculated with the help of the
identities $(2.10.1)$ and $(2.10.2)$. Suppose that $a$ and $y$
belong to $D$, $s=a^{\psi }$, $b=a^{\tau }$, $s_2=y^{\psi }$,
$b_2=y^{\tau }$. Then $(a^{\tau }y)=b(s_2b_2)$. According to
$(4.3.1)$ there exists a unique decomposition $b(s_2b_2)=s_3b_3$,
where $s_3\in A$, $b_3\in V$, hence $(a^{\tau }y)^{\tau }=b_3$. On
the other hand, by $(2.1.9)$
$ay=s(b(s_2b_2))t(s,b,y)=(ss_3)b_3t(s,b,y)/t(s,s_3,b_3)$. We denote
a subgroup $Z(D)\cap A$ in $Z(D)$ by $Z_A(D)$ or shortly $Z_A$, when
$D$ is specified. From Lemma 2.3 and $(4.1.4)$ it follows that
\par $(4.4.4)$ $Z(D)^{\tau }$ is isomorphic with $Z(D)/Z_A$, \par where $Z(D)^{\tau
}= \{ a^{\tau }: ~ a\in Z(D) \} $.
\par Let $Z_m(A)$ be a minimal subgroup in $Z(A)$ generated by a set
$ \{ t_A(a,b,c): ~ a\in A, ~b\in A, ~ c\in A \} $. From $(2.1.9)$ it
follows that $Z_m(A)\subset Z_A(D)$ and $AZ(D)$ is a submetagroup in
$D$. By virtue of Theorem 2.4 $(AZ(D))/Z_A(D)$ and $A/Z_A(D)$ are
groups such that $A/Z_A(D)\hookrightarrow (AZ(D))/Z_A(D)$. For each
$d\in D$ there exists a unique decomposition
\par $(4.4.5)$ $d=d^{\psi }d^{\tau }$ \\ by $(4.4.1)$. Take in
particular $\gamma \in Z(D)$, then $\gamma = \gamma ^{\psi }\gamma
^{\tau }$, where $\gamma ^{\psi }\in Z_A(D)$, $\gamma ^{\tau }\in
V$. Therefor $Z(D)/Z_A(D)\subset V$ and there exists a subset $V_0$
in $V$ such that $(Z(D)/Z_A(D))V_0=V$, since $Z(D)/Z_A(D)$ is a
subgroup in $(AZ(D))/Z_A(D)$ (see $(4.4.4)$). Formula $(4.4.5)$
implies that $(d^{\tau })^{\psi }=e$ and $(d^{\psi })^{\tau }=e$ for
each $d\in D$. Using this we subsequently deduce that
\par $(4.4.6)$ $(d^{\psi }\gamma )^{\psi }=d^{\psi }\gamma ^{\psi
}$,
\par $(4.4.7)$ $(d^{\psi }\gamma )^{\tau }=\gamma ^{\tau
}$,
\par $(4.4.8)$ $(d^{\tau }\gamma )^{\psi }=\gamma ^{\psi
}$,
\par $(4.4.9)$ $(d^{\tau }\gamma )^{\tau }=d^{\tau }\gamma ^{\tau
}$ \\ for each $d\in D$ and $\gamma \in Z(D)$. Hence \par $(d\gamma
)=(d\gamma )^{\psi } (d\gamma )^{\tau }=(d^{\psi }d^{\tau })(\gamma
^{\psi }\gamma ^{\tau })=(d^{\psi }\gamma ^{\psi }) (d^{\tau }\gamma
^{\tau })= (d^{\psi }\gamma )^{\psi } (d^{\tau }\gamma )^{\tau }$,
\\ where $d^{\psi }\gamma ^{\psi }\in A$ and $d^{\tau }\gamma ^{\tau
}\in V$. From a uniqueness of this representation it follows that
\par $(4.4.10)$ $(d\gamma )^{\psi }=d^{\psi }\gamma ^{\psi }$
and
\par $(4.4.11)$ $(d\gamma )^{\tau }=d^{\tau }\gamma ^{\tau }$
for each $d\in D$ and $\gamma \in Z(D)$.
\par Using $(4.4.11)$ we infer that
\par $(4.4.12)$ $(a^{\tau }y)^{\tau } = (ay)^{\tau }[t_D(a^{\psi
},(a^{\tau }y)^{\psi },(ay)^{\tau })/t_D(a^{\psi },a^{\tau
},y)]^{\tau }$. \par On the other hand, if $\gamma \in Z(D)$, then
$\gamma ^{\psi }=\gamma /\gamma ^{\tau }$ and Formulas $(4.4.12)$
and $(4.4.8)$ imply particularly that \par $(4.4.13)$ $(a^{\tau
}\gamma )^{\tau }=(a\gamma )^{\tau }$ for each $a\in D$ and $\gamma
\in Z(D)$, since $t_D(a,d,\gamma )=e$ for each $a$ and $d$ in $D$
and $\gamma \in Z(D)$. Then from $s=a^{\psi }$, $a^{\tau }y=s_3b_3$
it follows that $a^{\psi }(a^{\tau }y)^{\psi }=ss_3$ and $(ay)^{\psi
} =[(ss_3)b_3t_D(s,b,y)/t_D(s,s_3,b_3)]^{\psi }$, consequently, by
Lemma 2.3 and $(4.4.10)$
\par $(4.4.14)$ $a^{\psi }(a^{\tau }y)^{\psi }=(ay)^{\psi
}[t_D(a^{\psi },(a^{\tau }y)^{\psi },(ay)^{\tau })/t_D(a^{\psi
},a^{\tau },y)]^{\psi }$ \\ for each $a$ and $y$ in $D$.
Particularly
\par $(4.4.15)$ $a^{\psi }(a^{\tau }\gamma )^{\psi }=(a\gamma )^{\psi }$
for each $a\in D$ and $\gamma \in Z(D)$. \par From $(4.4.12)$ and
$(4.4.13)$ it follows that the metagroup $D$ acts on $V$
transitively by right shift operators $R_y$, where $R_ya=ay$ for
each $a$ and $y$ in $D$. Therefore, we put
\par $(4.4.16)$ $(a^{\tau })^{[c]}:=(a^{\tau }c)^{\tau }$ for each
$a $ and $c$ in $D$. \par Then from $(4.4.12)$, $(4.4.13)$,
$(4.4.16)$, $(2.1.9)$ and Lemma 2.3 we deduce that for each $a$,
$c$, $d$ in $D$
\par $(4.4.17)$ $(a^{\tau })^{[cd]} =((a^{\tau
})^{[c]})^{[d]}[t_D((a^{\tau }c)^{\psi },(a^{\tau }c)^{\tau
},d)/$\par $(t_D((a^{\tau }c)^{\psi },((a^{\tau }c)^{\tau }d)^{\psi
},((a^{\tau }c)d)^{\tau }) t_D(a^{\tau },c,d))]^{\tau }$. \par In
particular, $(a^{\tau })^{[e]}=a^{\tau }$ for each $a\in D$. Next we
put $e^{\tau }=b_*$. It is convenient to choose $b_*=e$. Hence
$b_*^{[s]}=(e^{\tau })^{[s]}=(e^{\tau }s)^{\tau }=s^{\tau
}=e=e^{\tau }$ for each $s\in A$. Thus the submetagroup $A$ is the
stabilizer of $e$ and $(4.4.16)$ implies that \par $(4.4.18)$
$e^{[s]}=e$ and $e^{[q]}=q$ for each $s\in A$ and $q\in V$.

\par {\bf 4.5. Remark.} Let $B$ and $D$ be metagroups, $A$ be a
submetagroup in $D$, $V$ be a right transversal of $A$ in $D$. Let
also Conditions $(3.2.1)$-$(3.2.8)$ be satisfied for $A$ and $B$. By
Theorem 3.1 there exists a metagroup
\par $(4.5.1)$ $F=B^V$, where $B^V=\prod_{v\in V}B_v$, $~B_v=B$ for
each $v\in V$. \par It contains a submetagroup \par $F^*=\{ f\in F:
~ card (\sigma (f))<\aleph _0 \} $, \\ where $\sigma (f)= \{ v\in V:
~ f(v)\ne e \} $ is a support of $f\in F$, $card (\Omega )$ denotes
the cardinality of a set $\Omega $.
\par Let $T_hf=f^h$ for each $f\in F$ and $h: V\to A$.
We put \par $\hat{S}_d(T_hfJ)=T_{hS_d^{-1}} fS_dJ$, \\ where $J:
V\times F\to B$, $ ~ J(f,v)=fJv$, $~ S_dJv=Jv^{[d\setminus e]}$ for
each $d\in D$, $f\in F$ and $v\in V$. Then for each $f\in F$, $d\in
D$ we put
\par $(4.5.2)$ $f^{ \{ d \} } =\hat{S}_d(T_{g_d}fE)$, \\
where \par $(4.5.3)$ $s(d,v)=e/(v/d)^{\psi }$, $~g_d(v)=s(d,v)$, \\
$~fEv=f(v)$ for each $v\in V$ (see $(4.3.1)$ and $(4.4.16)$). Hence
\par $(4.5.4)$ $f^{ \{ e \} } =f$, \\ since $v^{e\setminus e}=v$ and
$s(e,v)=e$.

\par {\bf 4.6. Lemma.} {\it Let Conditions of Remark 4.5 be satisfied.
Then for each $f$ and $f_1$ in $F$, $d$ and $d_1$ in $D$, $v\in V$:
\par $(4.6.1)$ $(ff_1)^{ \{ d \} } (v) = \kappa (s(d,v),
f(v^{[d\setminus e]}),f_1(v^{[d\setminus e]})) f^{ \{ d \} }(v)
f_1^{ \{ d \} } (v)$ and \par $(4.6.2)$ $f^{ \{ dd_1 \} } (v)= \{
[(f^{\{ d_1 \} })^{ \{ d \} } ]^{w_2(d,d_1,v)} (v w_1(d,d_1,v))
\} w_3(d,d_1,v)$, \\
where $w_j=w_j(d,d_1,v)\in Z(D)$, $j\in \{ 1, 2, 3 \} $, $w_1^{\tau
}=w_1$.}
\par {\bf Proof.} Formulas $(4.5.2)$ and $(3.2.6)$ imply Identity
$(4.6.1)$. \par Let $v\in V$, $d$ and $d_1$ belong to $D$, $f\in F$,
then from $(4.5.2)$ and $(4.5.3)$ it follows that
\par $(4.6.3)$ $f^{ \{ dd_1 \} } (v)=f^{s(dd_1,v)}(v^{[(dd_1)\setminus
e]})$ and
\par $(4.6.4)$ $(f^{ \{ d_1 \} })^{ \{ d \}
}(v)= (f^{s(d_1,v)})^{s(d,v^{[d_1\setminus e]})}((v^{[d_1\setminus
e]})^{[d\setminus e]}) $.
\par From Formulas $(2.6.2)$, $(4.4.17)$, $(4.4.6)$, $(4.4.9)$, $(2.2.1)$,
$(2.2.3)$ and Lemma 2.3 we deduce that
\par $(dd_1)\setminus e = (d_1\setminus e)(d\setminus
e)t_D(dd_1,d_1\setminus e,d\setminus e)/t_D(d,d_1,d_1\setminus e)$
and
\par $(4.6.5)$ $v^{[(dd_1)\setminus e]}=(v^{[d_1\setminus
e]})^{[d\setminus e]}w_1(d,d_1,v)$, where
\par $w_1(d,d_1,v)=\gamma ^{\tau }$, where
\par $\gamma = t_D(dd_1,d_1\setminus e,d\setminus
e)t_D((v/d_1)^{\psi },(v/d_1)^{\tau },d\setminus e)$\par
$/[t_D(d,d_1,d_1\setminus e) t_D(v,e/d_1,e/d) t_D((v/d_1)^{\psi },
((v/d_1)^{\tau }/d)^{\psi }, (v/(dd_1))^{\tau })]$.
\par Then Formulas $(4.5.3)$, $(4.4.14)$, $(2.6.3)$,
$(2.2.3)$, $(4.4.10)$ and Lemma 2.3 imply that
\par $(4.6.6.)$ $s(dd_1,v)=e/([(v/d_1)(e/d)]^{\psi }\gamma _1^{\psi
}$, where
\par $\gamma _1=t_D(v/(dd_1),d,d_1)t_D(e/d,d,d\setminus
e)/t_D(v/(dd_1),d,d\setminus e)$ by $(2.2.3)$, $(2.6.3)$ and Lemma
2.3;
\par $(4.6.7)$ $[(v/d_1)(e/d)]^{\psi }=(v/d_1)^{\psi }[(v/d_1)^{\tau
}(e/d)]^{\psi }\gamma _2^{\psi }$, where
\par $\gamma _2=t_D((v/d_1)^{\psi },(v/d_1)^{\tau
},e/d)/t_D((v/d_1)^{\psi },((v/d_1)^{\tau }(e/d))^{\psi
},((v/d_1)(e/d))^{\tau })$.
\par Note that \par $(4.6.8)$ $s(dd_1,v)=(e/[(v/d_1)^{\tau
}(e/d)]^{\psi })(e/(v/d_1)^{\psi })\gamma _3/\{ \gamma _1^{\psi
}\gamma _2^{\psi } \} $ \\ by $(4.6.7)$ and $(2.6.1)$, where \par
$\gamma _3=t_D(e/(v/d_1)^{\psi },(v/d_1)^{\psi },[(v/d_1)^{\tau
}(e/d)]^{\psi })$\par $/t_D(e/[(v/d_1)^{\tau }(e/d)]^{\psi
},e/(v/d_1)^{\psi }, (v/d_1)^{\psi }[(v/d_1)^{\tau }(e/d)]^{\psi
})$. Then
\par $(v/d_1)^{\tau }=[v(d_1\setminus e)]^{\tau }\gamma _4^{\tau }=
v^{[d_1\setminus e]}\gamma _4^{\tau }$ \\ by $(2.2.1)$, $(2.2.3)$,
$(4.4.11)$ and $(4.4.16)$, where $\gamma _4\in Z_m(D)$. Hence
\par $(4.6.9)$ $[(v/d_1)^{\tau }(e/d)]^{\psi }=[v^{[d_1\setminus
e]}(e/d)]^{\psi }\gamma _5^{\psi }$, \\ since $(\gamma _4^{\tau
})^{\psi }=e$, where
\par $\gamma _5=t_D(v^{[d_1\setminus e]}/d,d,d\setminus
e)/t_D(e/d,d,d\setminus e)$.
\par Thus Identities $(4.6.6)$-$(4.6.9)$ imply that
\par $(4.6.10)$ $s(dd_1,v)=s(d,v^{[d_1\setminus e]}) s(d_1,v)
w_2(d,d_1,v)$, where
\par $w_2(d,d_1,v)=\gamma _3/(\gamma _1^{\psi }\gamma _2^{\psi
}\gamma _5^{\psi })$, $~w_2(d,d_1,v)\in Z(D)$. By Lemmas 2.2, 2.3
and Formula $(4.5.3)$ representations of $\gamma _j$ simplify:
\par $\gamma _2=t_D(e/s(d_1,v),v^{[d_1\setminus
e]},e/d)/t_D(e/s(d_1,v),e/s(d,v^{[d_1\setminus
e]}),v^{[(dd_1)\setminus e]})$,
\par $\gamma _3=t_D(s(d_1,v),e/s(d_1,v),e/s(d,v^{[d_1\setminus
e]}))$\par $/t_D(s(d,v^{[d_1\setminus
e]}),s(d_1,v),e/(s(d,v^{[d_1\setminus e]})s(d_1,v)))$.
\par Therefore $w_2(d,d_1,v)\in Z(D)\cap A$ for each $d$ and $d_1$ in $D$
and $v\in V$, since $s(d,d_1,v)$, $s(d,v^{[d_1\setminus e]})$,
$s(d_1,v)$ belong to $A$. Then from $(4.6.3)$, $(4.6.10)$,
$(3.2.4)$, $(3.2.5)$ we infer that
\par $(4.6.11)$ $f^{ \{ dd_1 \} }(v)=$\par $\{ [(f^{s(d_1,v)})^{s(d,v^{[d_1\setminus e]})}
]^{w_2(d,d_1,v)} ((v^{[d_1\setminus e]})^{[d\setminus e]}
w_1(d,d_1,v)) \} w_3(d,d_1,v)$, \\ where $w_3(d,d_1,v)=e/[\eta
(s(d,v^{[d_1\setminus e]})w_2,s(d_1,v),b)\eta (s(d,v^{[d_1\setminus
e]}),w_2,b^{s(d_1,v)})]$,
\par $w_2=w_2(d,d_1,v)$, $~b=f((v^{[d_1\setminus e]})^{[d\setminus
e]} w_1(d,d_1,v))$. Formulas $(4.4.16)$ and $(4.4.9)$ imply that
$(v\gamma )^{[a]}=v^{[a]}\gamma ^{\tau }$ for each $v\in V$ and
$\gamma \in Z(D)$, $a\in D$, consequently, $((vw_1)^{[d_1\setminus
e]})^{[d\setminus e]}=(v^{[d_1\setminus e]})^{[d\setminus e]}w_1$
and hence $vw_1\in V$ for each $v\in V$, $d$ and $d_1$ in $D$,
$w_1=w_1(d,d_1,v)$ by Formula $(4.6.5)$, since $w_1=\gamma ^{\tau
}$. Thus Formula $(4.6.2)$ follows from $(4.6.4)$ and $(4.6.11)$.

\par {\bf 4.7. Definition.}  Suppose that the conditions of Remark
4.5 are satisfied and on the Cartesian product $C=D\times F$ (or
$C^*=D\times F^*$) a binary operation is given by the following
formula:
\par $(4.7.1)$ $(d_1,f_1)(d,f)=(d_1d, \xi ((d_1^{\psi
},f_1),(d^{\psi },f)) f_1f^{ \{ d_1 \} })$, where
\par $\xi ((d_1^{\psi },f_1),(d^{\psi },f))(v)=\xi ((d_1^{\psi
},f_1(v)),(d^{\psi },f(v)))$ for every $d$ and $d_1$ in $D$, $f$ and
$f_1$ in $F$ (or $F^*$ respectively), $v\in V$.

\par {\bf 4.8. Theorem.} {\it Let $C$, $C^*$, $D$, $F$, $F^*$ be the
same as in Definition 4.7. Then $C$ and $C^*$ are loops and there
are natural embeddings $D\hookrightarrow C$, $F\hookrightarrow C$,
$D\hookrightarrow C^*$, $F^*\hookrightarrow C^*$ such that $F$ (or
$F^*$) is an almost normal subloop in $C$ (or $C^*$ respectively).}
\par {\bf Proof.} The operation $(4.7.3)$ is single-valued. Let
$a=(d,f)$ and $b=(d_0,f_0)$, where $d$ and $d_0$ are in $D$, $f$ and
$f_0$ are in $F$ (or $F^*$). \par The equation $ay=b$ is equivalent
to $dd_2=d_0$ and \\ $\xi ((d^{\psi },f),(d_2^{\psi },f_2))ff_2^{ \{
d \} }=f_0$, where $d_2\in D$, $f_2\in F$ (or $f_2\in F^*$
respectively), $y=(d_2,f_2)$, $~\xi ((d^{\psi },f),(d_2^{\psi
},f_2))(v)=\xi ((d^{\psi },f(v)),(d_2^{\psi },f_2(v)))$ for each
$v\in V$. Therefore $d_2=d\setminus d_0$, $~f_2^{ \{ d \} }=[\xi
((d^{\psi },f),((d\setminus d_0)^{\psi },f_2))f]\setminus f_0$ by
$(2.1.1)$ and Theorem 3.1. On the other hand, $f_2^{ \{ e \} }=f_2$
by $(4.5.4)$ and $f_2(v)=\{ [(f_2^{ \{ d \} })^{ \{ d_3 \} }]^{w_2}
(vw_1) \} w_3$ by $(4.6.2)$, where $w_j=w_j(d,d_3,v)$, $ ~ j\in \{
1, 2, 3 \} $, $~d_3=d\setminus e$ and $dd_3=e$ by $(2.2.4)$. Thus
using $(3.2.8)$ we get that
\par $y=(d\setminus d_0,\{ \{ ([\xi ((d^{\psi },f),((d\setminus d_0)^{\psi
},[(f\setminus f_0)^{ \{ d\setminus e \} }]^{w_2}w_3))f]\setminus
f_0)^{ \{ d\setminus e \} }\} ^{w_2}(vw_1) \} w_3)$ \\ belongs to
$C$ (or $C^*$ respectively) giving $(2.1.1)$.

\par Then we seek a
solution $x\in C$ (or $x\in C^*$ respectively) of the equation
$xa=b$. It is equivalent to two equations $d_1d=d_0$ and \\ $\xi
((d_1^{\psi },f_1(v)),(d,f(v)))f_1(v) f^{ \{ d_1 \} }(v)=f_0(v)$ for
each $v\in V$, where $d_1\in D$, $f_1\in F$ (or $f_1\in F^*$
respectively), $x=(d_1,f_1)$. Therefore, $d_1=d_0/d$ and
$f_1(v)=f_0(v)/[\xi (((d_0/d)^{\psi },f_1(v)),(d,f(v))) f^{ \{ d_0/d
\} }(v)]$. Thus \par $x=(d_0/d,f_0/ [\xi (((d_0/d)^{\psi
},f_1),(d,f)) f^{ \{ d_0/d \} }])$ \\ belongs to $C$ (or $C^*$
respectively) giving $(2.1.2)$. \par Moreover, $(e,e)(d,f)=(d,f)$
and $(d,f)(e,e)=(d,f)$ for each $d\in D$, $f\in F$ (or $f\in F^*$
respectively) by $(3.2.8)$ and $(4.7.3)$. Therefore Condition
$(2.1.3)$ is also satisfied. Thus $C$ and $C^*$ are loops.
\par Evidently $D\ni d\mapsto (d,e)$ and $F\ni f\mapsto (e,f)\in C$
(or $F^*\ni f\mapsto (e,f)\in C^*$ respectively) provide embeddings
of $D$ and $F$ (or $D$ and $F^*$ respectively) into $C$ (or $C^*$
respectively).
\par It remains to verify that $F$ (or $F^*$ respectively) is an
almost normal subloop in $C$ (or $C^*$ respectively). Assume that
$d_1\in D$, $f_1\in F$. Then \par $(d_1,f_1)F= \{ (d_1,\xi
((d_1^{\psi },f_1),(e,f))f_1f^{ \{ d_1 \} }): ~ f\in F \} $ and
\par $F(d_1,f_1)= \{ (d_1,\xi ((e,f),(d_1^{\psi }, f_1)) ff_1): ~ f\in F \} $. \\ Using the embedding $Z^V\hookrightarrow F$
and $(3.2.8)$ we infer that $(d_1,f_1)F=F(d_1,f_1)$, since $F^{ \{
d_1 \} }=F$ by $(4.4.16)$ and Lemma 4.1 and $(3.2.3)$. It can be
verified similarly that $F^*$ is the almost normal subloop in $C^*$.

\par {\bf 4.9. Definition.} The product $(4.7.1)$ in
the loop $C$ (or $C^*$) of Theorem 4.8 is called a smashed twisted
wreath product of $D$ and $F$ (or a restricted smashed twisted
wreath product of $D$ and $F^*$ respectively) with smashing factors
$\phi $, $\eta $, $\kappa $, $\xi $ and it will be denoted by
$C=D\Delta ^{\phi , \eta , \kappa , \xi }F$ (or $C^*=D\Delta ^{\phi
, \eta , \kappa , \xi }F^*$ respectively). The loop $C$ (or $C^*$)
is also called a smashed splitting extension of $F$ (or of $F^*$
respectively) by $D$.

\par {\bf 4.10. Theorem.} {\it Let the conditions of Remark 4.5 be satisfied
and $Z_m(D)\subseteq Z$, where $Z$ is as in $(3.2.1)$. Then $C$ and
$C^*$ supplied with the binary operation $(4.7.1)$ are metagroups.}
\par {\bf Proof.} In view of Theorem 4.8 $C$ and $C^*$ are
loops. To each element $b$ in $B$ there corresponds an element $\{
b(v): \forall v\in V, ~ b(v)=b \} $ in $F$ which can be denoted by
$b$ also. From Conditions $(3.2.2)$-$(3.2.8)$ we deduce that
\par $(4.10.1)$ $\gamma ^a=\gamma $ and $f^{\gamma }=f$ for every
$\gamma \in Z$ and $a\in A$. Hence $(4.10.1)$ and $(4.7.1)$ imply
that $(Z(A),Z(F))\subseteq Z(C)$. On the other hand, $w_1=\gamma
^{\tau }$ with $\gamma \in Z_m(D)$ and $w_2=\gamma _3/(\gamma
_1^{\psi }\gamma _2^{\psi }\gamma _5^{\psi })$ with $\gamma
_1$,...,$\gamma _5$ in $Z_m(D)$ (see $(4.6.10)$), hence the
condition $Z_m(D)\subset Z$ implies that Formula $(4.6.2)$
simplifies to
\par $(4.10.2)$ $f^{ \{ dd_1 \} }(v)=(f^{ \{ d \} }(v))^{ \{ d_1 \}
}w_3(d,d_1,v)$ \\ for each $f\in F$, $v\in V$, $d$ and $d_1$ in $D$,
since $Z\subseteq Z(A)$ by $(3.2.1)$. Next we consider the products
\par $(4.10.3)$ $I_1=((d_2,f_2)(d_1,f_1))(d,f)=((d_2d_1,\xi
((d_2^{\psi },f_2),(d_1^{\psi },f_1))f_2f_1^{ \{ d_2 \} })(d,f)$ and
\par $(4.10.4)$ $I_2=(d_2,f_2)((d_1,f_1)(d,f))=(d_2,f_2)(d_1d,\xi
((d_1^{\psi },f_1),(d^{\psi },f))f_1f^{ \{ d_1 \} })$.
\par Then $(4.7.1)$, $(4.10.4)$, $(3.2.6)$-$(3.2.8)$ imply that
\par $(4.10.5)$ $I_2=(d_2(d_1d),\xi ((d_1^{\psi },
f_1),(d^{\psi },f)) \xi ((d_2^{\psi },f_2), ((d_1d)^{\psi },f_1f^{
\{ d_1 \} }))$\par $ \kappa (s(d_2,v), f_1(v^{[d_2\setminus e]}),
f^{ \{ d_1 \} }(v^{[d_2\setminus e]}))f_2(v)[f_1^{ \{ d_2 \}
}(v)(f^{ \{ d_1 \} })^{ \{ d_2 \} }(v)]$.
\par From $(4.10.2)$ and $(4.10.3)$, $(4.6.2)$, $(3.2.8)$ we infer that
\par $(4.10.6)$ $I_1=((d_2d_1)d,\xi
((d_1^{\psi },f_1),(d^{\psi },f))\xi (((d_2d_1)^{\psi },f_2f_1^{ \{
d_2 \} }), (d^{\psi },f))$\par $ (f_2f_1^{ \{ d_2 \} })(f^{ \{ d_1
\} })^{ \{ d_2 \} }w_3$, \\ where $w_3=w_3(d_1,d_2,v)$. Therefore
from $(4.10.5)$ and $(4.10.6)$ we infer that
\par $(4.10.7)$ $I_1=t_C((d_2,f_2),(d_1,f_1),(d,f))I_2$, where
\par $t_C((d_2,f_2),(d_1,f_1),(d,f))=t_D(d_2,d_1,d)t_B(f_2,f_1,f^{ \{ d_2d_1 \}
}) \xi ((d_1^{\psi },f_1),(d^{\psi },f))$\par $ \xi ((d_2^{\psi
},f_2),((d_1d)^{\psi },f_1f^{ \{ d_1 \} }))\kappa
(s(d_2,v),f_1(v^{[d_2\setminus e]}),f^{ \{ d_1 \} }(v^{[d_2\setminus
e]})) $\par $/[\xi ((d_2^{\psi },f_2),(d_1^{\psi },f_1))\xi
(((d_2d_1)^{\psi },f_2f_1^{ \{ d_2 \} }), (d^{\psi
},f))w_3(d_1,d_2,v)]$;
\par $t_B(f_2,f_1,f)(v)=t_B(f_2(v),f_1(v),f(v))$;
\par $\xi ((d_2^{\psi },f_2),(d_1^{\psi },f_1))(v)=
\xi ((d_2^{\psi },f_2(v)),(d_1^{\psi },f_1(v)))$ for every $f$,
$f_1$, $f_2$ in $F$, $d$, $d_1$, $d_2$ in $D$, $v\in V$. Then from
Formula $(4.10.7)$, $Z(F)=(Z(B))^V$ (see Theorem 3.1) and $(3.2.1)$
it follows that the loops $C$ and $C^*$ satisfy Condition $(2.1.9)$,
since $(Z,Z^V)\subseteq Z(C)$. Thus $C$ and $C^*$ are metagroups.

\par {\bf 4.11. Note.} Generally if $A\ne \{ e \} $ and $A\ne D$,
$B$, $\phi $, $\eta $, $\kappa $, $\xi $ are nontrivial, where $A$,
$B$, $D$ are metagroups or particularly may be groups, then the
loops $C$ and $C^*$ of Theorem 4.8 can be non metagroups. If in
$(3.2.8)$ drop the conditions $\xi ((e,e),(v,b))=e$ and $\xi
((v,b),(e,e))=e$ for each $v\in V$ and $b\in B$, then the proofs of
Theorems 3.3, 3.4, 4.8 demonstrate that $C_1$ and $C_2$ are strict
quasi-groups, $C$ and $C^*$ are quasi-groups.

\par {\bf 4.12. Definition.} Let $P_1$ and $P_2$ be two loops
with centers $Z(P_1)$ and $Z(P_2)$. Let also
\par $(4.12.1)$ $\mu (a,b)=\nu (a,b) \mu (a) \mu (b)$ \\ for each $a$
and $b$ in $P_1$, where $\nu (a,b)\in Z(P_2)$. Then $\mu $ will be
called a metamorphism of $P_1$ into $P_2$. If in addition $\mu $ is
surjective and bijective, then it will be called a metaisomorphism
and it will be said that $P_1$ is metaisomorphic to $P_2$.

\par {\bf 4.13. Theorem.} {\it Suppose that $A$, $B$, $D$ are
metagroups, $A\subset D$, $V_1$ and $V_2$ are right transversals of
$A$ in $D$, $F_j=B^{V_j}$, \par $P_j=D\Delta ^{\phi , \eta , \kappa
, \xi }F_j$, $ ~ P_j^*=D\Delta ^{\phi , \eta , \kappa , \xi }F_j^*$,
$ ~ j\in \{ 1, 2 \} $. \\ Then $P_1$ is metaisomorphic to $P_2$ and
$P_1^*$ to $P_2^*$.}
\par {\bf Proof.} By virtue of Theorem 4.8 $P_j$ and $P_j^*$ are
loops, where $j\in \{ 1, 2 \} $, $~Z^{V_j}\subset Z(P_j)$. From
$(4.4.10)$ and $(4.5.3)$ it follows that
\par $(4.13.1)$ $s_j(\delta d,v)=s_j(d,v/\delta )=
\delta ^{\psi _j}s_j(d,v)$ \\ for each $d\in D$, $v\in V_j$ and
$\delta \in Z(D)$, where $s_j$, $v^{[a]_j}$, $d^{\tau _j}$, $d^{\psi
_j}$ correspond to $V_j$, $j\in \{ 1, 2 \} $. Then Formulas
$(4.4.16)$ and $(4.4.11)$ imply that \par $(4.13.2)$ $v^{[\delta
/d]_j}=v^{[e/d]} \delta ^{\tau _j}$ \\ for each $d\in D$, $v\in V_j$
and $\delta \in Z(D)$, $j\in \{ 1, 2 \} $. Therefore from Identities
$(4.13.1)$, $(4.13.2)$, $(4.6.10)$ and Lemma 2.3 we infer that
\par $(4.13.3)$ $w_2(\delta d, \delta _1 d_1, \delta
_2v)=w_2(d,d_1,v)$ \\ for each $d$ and $d_1$ in $D$, $\delta $,
$\delta _1$ and $\delta _2$ in $Z(D)$, $v\in V$.
\par For each $f\in F_1$ and $v\in V_2$ we put
\par $(4.13.4)$ $\mu f(v)=f^{e/v^{\psi _1}}(v^{\tau _1})$.
\par From Lemma 4.1 it follows that $V_2^{\tau _1}=V_1$ and
$v_1^{\tau _1}\ne v_2^{\tau _1}$ for each $v_1\ne v_2$ in $V_2$,
where $V_2^{\tau _1} = \{ v^{\tau _1}: v\in V_2 \} $. Then
$(4.10.1)$, $(4.4.10)$, $(4.13.3)$, $(4.13.4)$ and Lemma 2.2 imply
that \par $(4.13.5)$ $f^{ \{ d \} \mu }=f^{ \mu \{ d \} }$ \\
for each $f\in F_1$, $d\in D$, where $f^{ \mu \{ d \} }=(\mu f)^{ \{
d \} }$, $f^{ \{ d \} \mu }=\mu (f^{ \{ d \} })$ (see also
$(4.5.2)$). From Identity $(4.13.5)$ and Conditions $(3.2.6)$,
$(3.2.7)$ we infer that
\par $\mu ((d_1,f_1)(d,f))(v)=\kappa (e/v^{\psi _1},f_1(v^{\tau
_1}),f^{ \{ d_1 \} }(v^{\tau _1}))(\mu (d_1,f_1))(\mu (d,f)) $ \\
for each $d$ and $d_1$ in $D$, $f$ and $f_1$ in $F_1$, $v\in V_2$,
where $\mu (d,f)=(d,\mu f)$, $(d,f)(v)=(d,f(v))$. Hence \par $\nu
((d_1,f_1), (d,f))(v)=\kappa (e/v^{\psi _1},f_1(v^{\tau _1}),f^{ \{
d_1 \} }(v^{\tau _1}))\in Z$ \\ for each $v\in V_2$ (see also
$(3.2.1)$, $(3.2.3)$). Thus $P_1$ is metaisomorphic to $P_2$ and
$P_1^*$ to $P_2^*$.

\par {\bf 4.14. Theorem.} {\it Suppose that $D$ is a non-trivial metagroup.
Then there exists a smashed splitting extension $C^*$ of a
non-trivial central metagroup $H$ by $D$ such that $[H,C^*]Z(H)=H$,
where $[a,b]=(e/a)((e/b)(ab))$ for each $a$ and $b$ in $C^*$.}
\par {\bf Proof.} Let $d_0$ be an arbitrary fixed element in $D-Z(D)$.
Assume that $A$ is a submetagroup in $D$ such that $A$ is generated
by $d_0$ and a subgroup $Z_0$ contained in a center $Z(D)$ of $D$,
$~Z_m(D)\subseteq Z_0\subseteq Z(D)$, where $Z_m(D)$ is a minimal
subgroup in a center $Z(D)$ of $D$ such that $t_D(a,b,c)\in Z_m(D)$
for each $a$, $b$, $c$ in $D$. Therefore
\par $(4.14.1)$ $a^ka^n=p(k,n,a)a^{k+n}$ \\ for each $a\in A$,
$k$ and $n$ in ${\bf Z}= \{ 0, -1, 1, -2, 2,... \} $, where the
following notation is used $a^2=aa$, $a^{n+1}=a^na$ and
$a^{-n}=e/a^n$, $a^0=e$ for each $n\in {\bf N}$, $p(k,n,a)\in
Z_m(A)$. Hence in particular $A$ is a central metagroup. Then
$d_0Z_m(A)$ is a cyclic element in the quotient group $A/Z_m(A)$
(see Theorem 2.4). Then we choose a central metagroup $B$ generated
by an element $b_0$ and a commutative group $Z_1$ such that
$b_0\notin Z_1$, $ ~ Z_m(D)\hookrightarrow Z_1$ and
$Z(A)\hookrightarrow Z_1$, the quotient group $B/Z_m(B)$ is of
finite order $l>1$. Then let $\phi : A \to {\cal A}(B)$ satisfy
Condition $(3.2.3)$ and be such that \par $(4.14.2)$ $\phi
(d_0)b_0=b_0^2$.
\par To satisfy Condition $(4.14.2)$ a natural number $l$
can be chosen as a divisor of $2^{|d_0Z_m(A)|}-1$ if the order
$|d_0Z_m(A)|$ of $d_0Z_m(A)$ in $A/Z_m(A)$ is positive; otherwise
$l$ can be taken as any fixed odd number $l>1$ if $A/Z_m(A)$ is
infinite. \par Then we take a right transversal $V$ of $A$ in $D$ so
that $A$ is represented in $V$ by $e$. Let $\Xi $, $\eta $, $\kappa
$, $\xi $ be chosen satisfying Conditions $(3.2.2)$-$(3.2.8)$, where
$Z_m(B)\hookrightarrow Z$, $Z_m(A)\hookrightarrow Z$,
$Z_0\hookrightarrow Z$, $Z_1\hookrightarrow Z$. With these data
according to Theorem 4.10 $C^*$ is a metagroup, since
$Z_m(D)\hookrightarrow Z_1$ and $Z_m(D)\hookrightarrow Z_0$. That is
$C^*$ a smashed splitting extension of the central metagroup $F^*$
by $D$. \par Apparently there exists $f_0\in F^*$ such that
$f_0(e)=b_0$, $f_0(v)=e$ for each $v\in V - \{ e \} $. Therefore
$f_0^{ \{ v \} }(v)=b_0$ for each $v\in V$, since $s(v,v)=e$,
$v^{[v\setminus e]}=[v(v\setminus e)]^{\tau }=e$. \par Let $v_1\ne
v_2$ belong to $V$. Then $(v_2(v_1\setminus e))^{\tau }=v_3\in V$.
Assume that $v_3=e$. The latter is equivalent to $v_2(v_1\setminus
e)=a\in A$. From $(2.2.3)$ it follows that $v_2=a/(v_1\setminus
e)=\gamma av_1$, where $\gamma =t_D(v_1,v_1\setminus
e,v_1)/t_D(av_1,v_1\setminus e, v_1)$ by $(2.2.1)$ and Lemma 2.3,
since $ e/(v_1\setminus e)=v_1$. Hence $v_2=v_2^{\tau }=(\gamma
av_1)^{\tau }=\gamma ^{\tau }v_1$ by $(4.4.11)$ and consequently,
$(v_2(v_1\setminus e))^{\tau }=\gamma ^{\tau }=e$ contradicting the
supposition $v_1\ne v_2$. Thus $v_3\ne e$ and consequently, $f_0^{
\{ v_1 \} }(v_2)=e^{s(v_1,v_2)}=e$ by $(3.2.4)$. This implies that
$\{ f_0^{ \{ v \} }: v\in V \} Z(F^*)$ generate $F^*$.
\par Evidently $[v(d_0\setminus e)]^{\tau }\ne e$ for each $v\in V-
\{ e \} $, since $d_0\setminus e\in A$ and the following conditions
$s\in D$, $sq\in A$, $q\in A$ imply that $s\in A$, because $A$ is
the submetagroup in $D$. Note that $e/d=(d\setminus
e)/t_A(e/d,d,d\setminus e)$ for each $d\in A$ by $(2.2.1)$,
consequently, $s(d,e)=dt_A(e/d,d,d\setminus e)$. On the other hand
$t_A(a,b,c)\in Z$ for each $a$, $b$, $c$ in $A$ and
\par $(4.14.3)$ $f_0^{\gamma }=f_0$ for each $\gamma \in Z$ \\ by
$(4.10.1)$, hence $f_0^{ \{ d_0 \} }(e)=\phi (d_0)b_0=b_0^2$ and
consequently,
\par $(4.14.4)$ $f_0^{ \{ d_0 \} }=f_0^2$, \\ since $f_0^{ \{ d_0 \}
}(v)=e$ for each $v\in V- \{ e \} $.
\par Therefore we deduce using $(4.14.3)$ that
\par $[(e,f_0),(e/d_0,e)] = (e,wf_0)$, where \\ $w=\xi
((e,f_0),(e/d_0,e)) \xi ((d_0,e),(e/d_0,f_0)) \xi
((e,e/f_0),(e,(f_0)^2)) /t_{F^*}(e/f_0,f_0,f_0)$, \\
$~t_{F^*}(f,g,h)(v)=t_B(f(v),g(v),h(v))$ for each $v\in V$, $f$, $g$
and $h$ in $F^*$. Thus $w=w(v)\in Z$ for each $v\in V$ and $f_0\in
[F^*,C^*]$, since $Z^V\cap F^*\subset Z(F^*)$. Hence $F^*\subseteq
[F^*,C^*]Z(F^*)$, since $F^*\hookrightarrow C^*$ and $Z(C^*)\cap
F^*\subseteq Z(F^*)$. On the other hand, $Z_m(A)\hookrightarrow Z$,
$Z_m(B)\hookrightarrow Z$, $Z_m(D)\hookrightarrow Z_j$,
$Z_j\hookrightarrow Z$ for each $j\in \{ 0, 1 \} $. Therefore
$(4.14.3)$, $(4.14.4)$ and $(4.10.2)$ imply that $cF^*=F^*c$ and
$c[F^*,C^*]Z(F^*)=[F^*,C^*]Z(F^*)c$ for each $c\in C^*$. Hence
$[F^*,C^*]Z(F^*)\subseteq F^*$. Taking $H=F^*$ we get the assertion
of this theorem.

\par {\bf 4.15. Corollary.} {\it Let the conditions of Theorem 4.14
be satisfied and $D$ be generated by $Z_m(D)$ and at least two
elements $d_1$, $d_2$,... such that $d_1\ne e$ and $[d_2\setminus e,
d_1\setminus e]=e$. Then the smashed splitting extension $C^*$ can
be generated by $Z(F^*)$ and elements $c_1$, $c_2$,... such that
$d_j\setminus e \in F^*c_j$ for each $j$.}
\par {\bf Proof.} We take $d_0=d_1$ in the proof of Theorem 4.14,
$c_1=(d_1\setminus e,e)$, $c_2=(d_2\setminus e,f_0)$,
$c_j=(d_j\setminus e,e)$ for each $j\ge 3$. Therefore $(4.4.14)$,
$(4.14.4)$ and $(3.2.8)$ imply that
\par $[c_2,c_1]=(e, p f_0)$, where \par $p= \xi ((d_2\setminus e,f_0),
(d_1\setminus e,e))\xi ((d_1,e),((d_2\setminus e)(d_1\setminus
e),f_0)) $\par $\xi ((e,e)/(d_2\setminus e,f_0),(d_2\setminus
e,(f_0)^2)) /t_{F^*}(e/f_0,f_0,f_0)$, \\ since $[d_2\setminus e,
d_1\setminus e]=e$ and $e/(d_2\setminus e)=d_2$. Thus the
submetagroup of $C^*$ which is generated by $Z_m(D)$ and $ \{ c_j: j
\} $ contains the metagroup $D$ and $(e,pf_0)$. Therefore the
following set $\{ f^{ \{ d \} }: d\in D \} Z(F^*)$ generate the
central metagroup $F^*$, since $V\subset D$ and $\{ f^{ \{ v \} }:
v\in V \} Z(F^*)$ generate $F^*$. Notice that $Z_m(D)\hookrightarrow
Z(F^*)$. Hence $ \{ c_j: j \} Z(F^*)$ generate $C^*$.

\par {\bf 4.16. Conclusion.} The results of this article
can be used for further studies of metagroups, quasi-groups, loops
and noncommutative manifolds related with them. Besides applications
of metagroups, loops and quasi-groups outlined in the introduction
it is interesting to mention possible applications in mathematical
coding theory and its technical applications
\cite{blautrctb,petbagsychrtj,srwseabm14}, because frequently codes
are based on binary systems.

\end{document}